\documentclass[12pt,letterpaper]{amsart}

\usepackage{amssymb,latexsym,amsmath,amsthm,graphicx}
\numberwithin{equation}{section}

\newtheorem*{remark}{Remark}

\usepackage{xcolor}

\newcommand{\pd}{\partial}

\newcommand{\ol}{\overline}

\input epsf

\begin{document}

\title{A Lower Bound on the Growth of Minimal Graphs} 

\author[Allen Weitsman]{Allen Weitsman}

\address{Department of Mathematics\\ Purdue University\\
West Lafayette, IN 47907-1395}
\address{Email: weitsman@purdue.edu}


\begin{abstract}
We show that for minimal graphs in $R^3$ having 0 boundary values over simpy connected domains,
the maximum over circles of radius r must be at least of the order $r^{1/2}$. 


\end{abstract}

\maketitle

\section{Introduction} Let $\ D\ $ be an unbounded domain in $\Bbb R^2$.  
 We are interested in
solutions to the minimal surface equation  for the boundary value problem
\begin{equation}
\label{mineq1}
Lu=\text{div}\left({\nabla u\over\sqrt{1+|\nabla u|^2}}\right)=0
\end{equation}
in $\ D\ $ with
\begin{equation}
\label{mineq2}
u>0\quad\text{ in }\ \ D,\qquad  u=0\quad \text{ on }\  \partial D.
\end{equation} 

We shall use complex notation $z=x+iy$  for convenience.  With $M(r)$ being the maximum 
value of $u(z)$ on $D\cap\{|z|=r\}$, we have previously studied upper bounds on the growth rate of $M(r)$ 
under various conditions  \cite{LW}, \cite{W2}, \cite {W3}. We have also obtained some information on the lower bounds.

In \cite{LW} we observed the general result

{\bf Theorem A.}  \emph{Suppose $D$ is a domain with $\pd D \neq \emptyset$,
and $u$ as in (\ref{mineq1}) and (\ref{mineq2}).
Then $u(z)$ has at least logarithmic growth.}

From  \cite{W3} we have

{\bf Theorem B.} \emph{ Let $u$ satisfy (\ref{mineq1}) and (\ref{mineq2}) with $D$ simply 
connected and contained in a half plane.  Then,}
$$
\underset {r\to\infty}{\liminf}\,\frac{M(r)}{ r}>0.
$$
In this note we prove the following

{\bf Theorem 1.}  \emph{Suppose that $u(z)$ satisfies (\ref{mineq1}) and (\ref{mineq2}) with $D$  simply
connected.  Then}
\begin{equation}
\label{mingrowth}
\underset {r\to\infty}{\liminf}\,\frac{\log M(r)}{\log r}\geq 1/2.
\end{equation}

The catenoid \cite[p.161]{D} shows the necessity of simple connectivity in Theorem 1.

We note that (\ref{mingrowth}) with $\limsup$ in place of $ \liminf$ follows from \cite[Theorem 1]{W}).  Also, 
(\ref{mingrowth}) with $1/2$ replaced by $1/\pi$ can be deduced from the work of Miklyukov \cite[p.64]{M}.

The example given in \cite[p. 1085]{W2} shows that (\ref{mingrowth}) is sharp. (See also\cite[p 3391]{LW}
for related examples.)

\section{preliminaries}
Let $u$ be a solution to (\ref{mineq1}) and (\ref{mineq2}) over a 
simply connected domain $D$. 
We shall  make use of the parametrization of a surface given by $u$ in isothermal coordinates using 
Weierstrass functions $\left( x(\zeta), y (\zeta ), U (\zeta) \right)$
with $\zeta$ in the right half plane $H$.   Our notation will then be
given by 
\begin{equation}
\label{downstairs}
f(\zeta) = x(\zeta) +iy(\zeta)\quad\zeta = \sigma+i\tau\in H.
\end{equation}
Then $f(\zeta )$ is univalent and harmonic, and since $D$ is simply connected it can be written in
the form 
\begin{equation}
\label{decomp1}
f(\zeta) = h(\zeta) + \ol{g(\zeta)}
\end{equation}
where $h(\zeta )$ and $g(\zeta)$ are analytic in $H$, 
\begin{equation}
\label{Dilatation}
|h'(\zeta)|>|g'(\zeta)|.
\end{equation}

We dismiss the trivial case $g'\equiv 0$ and may assume for  later convenience that $f(0)=0$.
Regarding the height function, we have (cf. \cite[\S 10.2]{D})
\begin{equation}
\label{decomp2}
U(\zeta )= \pm\  2\Re e \, i\int \sqrt{h'(\zeta )g'(\zeta )}\,d\zeta .
\end{equation}

Now, $z=f(\zeta),\  u(f(\zeta))=U(\zeta)$ and $U(\zeta )$ is harmonic and positive in
 $H$ and vanishes on $\partial H$.  Thus,
(cf. \cite[p. 151]{T}), 
\begin{equation}
\label{height}
U(\zeta )= K\,\Re e\, \zeta,
\end{equation}
where $K$ is a positive constant.  This with
(\ref{decomp2}) gives
 
$$
 g'(\zeta) = - \frac{C}{h'(\zeta)}
$$
where $C$ is a positive constant.  By reparametrizing we may assume that
\begin{equation}
\label{rescale}
U(\zeta ) =2 \Re e\,\zeta \ \ \textrm{and}\ \   g'(\zeta )=-1/h'(\zeta ),
\end{equation}
and then the \emph{analytic dilatation} \cite[p.6]{D} $a(\zeta)$ satisfies
\begin {equation}
\label{dilatation1}
a(\zeta)=-1/h'(\zeta)^2.
\end{equation}
Furthermore,  from (\ref{Dilatation}) we have, in particular, that
\begin{equation}\label{decomp3}
|h'(\zeta)|=1/|g'(\zeta)|>1.
\end{equation}

The strategy  will be to analyze   $f(\zeta)$ in  sectors 
$$
S_\epsilon =\{(-\pi +\varepsilon)/2 <\arg \zeta<(\pi -\varepsilon)/2\},
$$
 where $0< \varepsilon <\pi/2$. 
 We also define, for  fixed $\rho>1$, 
\begin{equation}
\label{rectangle1}
 S_{\varepsilon,n}=S_\varepsilon\cap\{\rho^n\leq |\zeta | \leq \rho^{n+1}\}\ n=0,1,2,....
\end{equation}

\section{Quasiconformal mappings}
We shall have occasion to view the harmonic mapping  described in \S 2 as a  quasiconformal mapping.  A one to one sense preserving mapping  $f$ in a domain $D$ is quasiconformal, if its
 \emph {complex dilatation} $ \delta (\zeta )$ defined by  (cf.(\cite[p. 5]{D}))
\begin{equation}
\label{dilatation}
\delta(\zeta) = \frac{f_{\overline{\zeta}}(\zeta)}{f_\zeta(\zeta)},
\end{equation}
satisfies
\begin{equation}
\label{quasiconformal}
\underset{\zeta \in D}{\sup}\,|\delta(\zeta)|<1.
\end{equation}
Henceforth, we shall refer to $|\delta(\zeta)|$ simply as the dilatation.

The dilatation is a conformal invariant, and the inverse mapping has the same dilatation at corresponding
points \cite[p. 9]{A}.

We shall need a modification of the Ahlfors distortion theorem which requires slight changes in the
standard proof \cite[pp. 94-97]{Ne}. 

In the classical setting we have a simply connected region $G$ with accessible boundary
points $Z_1=X_1+iY_1$ and $Z_2=X_2+iY_2$ 
$Z_1$ and $Z_2$.   We assume that $-\infty\leq X_1=\inf \Re e\,z$ for $z\in G$ 
and $\infty\geq X_2=\sup \Re e\, z$ for $z\in G$. We consider 
 $z=x+iy$ in $G$ with 
cross cuts $\Theta_x$  separating $Z_1$ and $Z_2$  in $G$ (See (\cite[pp. 94-95]{Ne}) for more details).  Let $\Theta
(x)$ be the length of $\Theta_x$. Let $w(z)=\mu(z)+i\nu(z)$ be a conformal mapping of $G$ onto the strip
$\{|\nu|<a/2\}$  such that $Z_1$ corresponds to $-\infty$ and $Z_2$ to $+\infty$.

 If $\mu_1(x)$ denotes the 
smallest value on the cross cut and $\mu_2(x)$ the largest, then
the classical distortion theorem is as follows.

{\bf Theorem  A. } \emph{If
$$
\int_{x_1}^{x_2} \frac{dx}{\Theta(x)}>2,
$$
then 
$$
\mu_1(x_2)-\mu_2(x_1)\geq a\int_{x_1}^{x_2} \frac {dx} {\Theta(x)}-4a.
$$
}

For our purposes, the strip $\Sigma_\varepsilon$ will be the (principal branch) logarithmic image of $S_\varepsilon$ in the $w=\mu+i\nu$ plane and $G$ will be  
 the image of a fixed branch of $\log f(S_\varepsilon)$ in the $z=x+iy$ plane  
with $f$ as in \S2.  As previously mentioned, we assume for convenience that $f(0)=0$ so that in $G$, $\Re e\, z$  extends from $-\infty$ to $+\infty$.

Let $w(z)=\log(f^{-1}(e^z))$ for the principal branch of log which then has the same dilatation as $f$
  at corresponding points.

{\bf Lemma 1.}  \emph{With the above notations,  
let $R$
be a rectangle in the $\mu +i\nu$ plane
$$
R=\Sigma_\varepsilon\cap\ \{\alpha \leq \mu\leq \beta\}\quad 0<\alpha<\beta,
$$
and suppose that  $w(z)$ has dilatation less than $\delta_0$ in $w^{-1}(R)$. 
Then for $x_1+iy_1$ and $x_2+iy_2$ in $w^{-1}(R)$ ($x_1<x_2$) and $a=\pi-\varepsilon$, if
\begin{equation}
\label{Ahlfors1}
 \int_{x_1}^{x_2} \frac{dx}{\Theta(x)}>2\frac{1+\delta_0}{1-\delta_0}, 
\end{equation}
}
then
\begin{equation}
\label{Ahlfors4}
\mu_1(x_2)-\mu_2(x_1)\geq a\frac{(1-\delta_0)}{(1+\delta_0)}\int_{x_1}^{x_2}\frac{dx}{\Theta(x)} 
-4a.
\end{equation}

{\bf Proof.  }  Our proof follows \cite[pp 95-97]{Ne}.
  
The length of the arc $L_x$ corresponding to $w(\Theta_x)$ is at least
 $\sqrt{a^2+\omega(x)^2}$, where  $\omega(x)=\mu_2(x)-\mu_1(x)$.   Also,
$$
L_x\leq \int_{\Theta_x}(|w_z|+|w_{\overline z}|)dy\leq\int_{\Theta_x}(|w_z|(1+\delta_0)dy
\leq \sqrt{\int_{\Theta_x}dy\int_{\Theta_x}|w_z|^2(1+\delta_0)^2dy}.
$$
  Thus,
  $$
  a^2+\omega(x)^2\leq \Theta(x)\int_{\Theta_x}|w_z|^2(1+\delta_0)^2dy.
  $$
Then,
$$
a^2\int_{x_1}^{x_2}\frac{dx}{\Theta(x)}dx+\int_{x_1}^{x_2}\frac{\omega(x)^2dx}{\Theta(x)}dx 
\leq\int_{x_1}^{x2}\int_{\Theta_x}|w_z|^2(1+\delta_0)^2dydx
$$
$$
=(1+\delta_0)^2\int_{x_1}^{x_2}\int_{\Theta_x}(|w_z|^2-|w_{\overline z}|^2+|w_{\overline z}|^2)dydx
 $$ 
 $$\leq (1+\delta_0)^2\int_{x_1}^{x_2}\int_{\Theta_x}(|w_z|^2-|w_{\overline z}|^2
 +\frac{\delta_0^2|w_z|^2}{(|w_z|^2-|w_{\overline z}|^2)}(|w_z|^2-|w_{\overline z}|^2)dydx
 $$
 $$
 \leq (1+\delta_0)^2\int_{x_1}^{x_2}\int_{\Theta_x}(|w_z|^2
 -|w_{\overline z}|^2)(1+\frac{\delta_0^2}{1-\delta_0^2})dydx
 $$
 $$
 =\frac {1+\delta_0} {1-\delta_0}\int_{x_1}^{x_2}\int_{\Theta_x}(|w_z|^2-|w_{\overline z}|^2)dydx
 $$
 $$
 \leq a\frac{(1+\delta_0)}{(1-\delta_0)}(\mu_2(x_2)-\mu_1(x_1))
$$
$$
=a\frac  {(1+\delta_0)} {((1-\delta_0)}(\mu_1(x_2)-\mu_2(x_1) +\omega(x_2)+\omega(x_1)).
$$
Summarizing this we have
\begin{equation}
\label{Ahlfors2}
\mu_1(x_2)-\mu_2(x_1)\geq \frac{1-\delta_0}{1+\delta_0}\left(a\int_{x_1}^{x_2}\frac{dx}{\Theta(x)}+\frac 1 a \int_{x_1}^{x_2}\frac{\omega(x)^2}{\Theta(x)}dx\right) -\omega(x_1)-\omega(x_2).
\end{equation}
For $x_0\in (x_1,x_2)$, let
$$
\lambda(x)=\frac 1 a\frac{(1-\delta_0)}{(1+\delta_0)}\int _{x_1}^x\frac{\omega (t)^2}{\Theta (t)}dt.
$$
Let $m>0$ be fixed and $\mathcal{ E}$ be the set of $x>x_0$ such that $\lambda(x)<\omega(x)-m$.
Then,
$$
a\frac{1+\delta_0)}{(1-\delta_0)}\Theta((x)\frac{d\lambda}{dx}=\omega(x)^2>(\lambda(x)+m)^2,
$$
so that
$$
\int_{\mathcal {E}}\frac{dx}{\Theta(x)}\leq a\frac{(1+\delta_0)}{(1-\delta_0)}\int_{\mathcal{E}}
\frac{d \lambda}{(\lambda+m)^2}\leq\frac a m\frac{(1+\delta_0)}{(1-\delta_0)}.
$$
Similarly, if $\mathcal{F}$ is the set of $x<x_0$ such that $-\lambda (x)<\omega(x)-m$,
$$
\int_{\mathcal{F}}\frac{dx}{\Theta(x)}\leq \frac a m\frac{(1+\delta_0)}{(1-\delta_0)}.
$$ 
Choose $x_0$ such that
$$
\int_{x_1}^{x_0}\frac{dx}{\Theta(x)} =\int_{x_0}^{x_2}\frac{dx}{\Theta(x)}= \frac 1 2 \int_{x_1}^{x_2}\frac{dx}{\Theta(x)}.
$$
Then, if 
$$\
{\int_{x_1}^{x_2}\frac{dx}{\Theta(x)}>2}\frac a m\frac{(1+\delta_0)}{(1-\delta_0)},
$$
 we may take
$x_1',\ x_2'\ (x_1<x_1'<x_2'<x_2)$ such that
\begin{equation}
\label{Ahlfors3}
\int_{x_1}^{x_1'}\frac{dx}{\Theta(x)}=\int_{x_2'}^{x_2}\frac{dx}{\Theta(x)}=\frac a m\frac{(1+\delta_0)}{(1-\delta_0)}.
\end{equation}
So, there exist $\xi_1\in (x_1,x_1')$ and $\xi_2\in (x_2',x_2)$ such that
$$
-\lambda(\xi_1)\geq\omega (\xi_1)-m\ \  \textrm{and}\  \ \lambda (\xi_2)\geq \omega (\xi_2)-m.
$$
Then
$$
\frac 1 a \frac{(1-\delta_0)}{(1+\delta_0)}\int_{\xi_1}^{\xi_2}\frac{\omega(x)^2}{\Theta(x)}dx
=\lambda(\xi_2)-\lambda(\xi_1)\geq\omega(\xi_1)+\omega(\xi_2)-2m.
$$
From (\ref{Ahlfors2})) we then have
$$
\mu_1(\xi_2)-\mu_2(\xi_1)\geq a\frac{(1-\delta_0)}{(1+\delta_0)}\int_{\xi_1}^{\xi_2}\frac{dx}{\Theta(x)}
+\omega(\xi_1)+\omega(\xi_2)-2m-\omega(\xi_1)-\omega(\xi_2).
$$
From this and (\ref{Ahlfors3}) we deduce that
$$
\mu_1(x_2)-\mu_2(x_1)\geq a\frac{(1-\delta_0)}{(1+\delta_0)}\int_{x_1}^{x_2}\frac{dx}{\Theta(x)} 
-\frac {2a^2} m-2m.
$$
With $m=a$ and then
$\displaystyle{\int_{x_1}^{x_2}\frac{dx}{\Theta(x)}>2\frac{1+\delta_0}{1-\delta_0}}$
we  have

$$
\mu_1(x_2)-\mu_2(x_1)\geq a\frac{(1-\delta_0)}{(1+\delta_0)}\int_{x_1}^{x_2}\frac{dx}{\Theta(x)} 
-4a.
$$
\qed
\section{The Parameters}

We now select parameters in order to utilize Lemma 1.

To begin with we fix $\varepsilon_1>0$ and take $a=\pi-\varepsilon_1$.
  Next, we
 fix $0<\delta_0<1/2$
so that $\displaystyle{\frac {1-\delta_0} {1+\delta_0}}>1-\varepsilon_1$ 
and for
\begin{equation}
\label{4.0}
C_1=\frac{2\pi}{a(1-\varepsilon_1)}
\end{equation}
define $\varepsilon_2$ by
\begin{equation}
\label{4.1}
\varepsilon_2=C_1-2.
\end{equation}
We then define
\begin{equation}
\label{4.2}
C_2=\exp(\frac{8\pi}{(1-\varepsilon_1)}
\end{equation}
 and fix a value $\rho$ in (\ref{rectangle1}) large enough so that 
\begin{equation}
\label{4.3}
\rho>e^{5\pi} \ \ \textrm{and}\ \ C_2/\rho^{\varepsilon_2}<1/2.
\end{equation}
The $R$ in Lemma 1 now corresponds to a sequence
$$
R_{\varepsilon_1n}=\Sigma_{\varepsilon_1}\cap \{(n-1)\log\rho\leq\mu\leq n\log\rho\}.
$$

Finally, we note that since $\log|h'|$ is a positive harmonic function in $H$, if $M$ is the maximum of
$\log| |h'|$ in $S_{\varepsilon_1,n}$ and $m$ the minimum of $\log |h'|$ in $S_{\varepsilon_1,n}$,  there exists a constant $K=
K(\varepsilon_1,\rho)$ (independent of $n$) such that $m/M\geq K$  Thus, we may fix a value
$M_0$  for $M$ such that if  $\log|h'|$ has maximum greater or equal to $M_0$, then the minimum of
$\log|h'|$ will be greater than $\log(1/\sqrt{\delta_0})$ and hence by (\ref{dilatation1})
$$
| \delta(\zeta)|=1/|h'(\zeta)|^2< \delta_0 \quad  \zeta\in S_{\varepsilon_1,n}.  
$$
The objective will be to bound $\log |f(\sigma)|$ for large $\sigma$, with $N$ chosen so that 
$$
\rho^{N-1}<\sigma<\rho^N.
$$

Note that $|f(\sigma)|$ is unbounded.  In fact since $U(\sigma)=2\sigma$ (recall (\ref{rescale})),
if $|f(\sigma)|$ were bounded, then (\ref{mingrowth}) would hold trivially.

\section{Proof of Theorem 1}
With the conventions in \S 4, we  distinguish two cases.

{\bf Case 1.}  Here we consider the case in which the maximum of $\log|h'|$ is at least $M_0$ in a
given $S_{\varepsilon_1n}$.  For the corresponding $R_{\varepsilon_1n}$,  then $w(z)=\log (f^{-1}(e^z))$ has dilatation less than $\delta_0$
in $w^{-1}(R_{\varepsilon_1n})$, and in Lemma 1 
we take $x_1+iy_1=w^{-1}((n-1)\log\rho)$ and $x_2+iy_2=w^{-1}(n\log\rho)$.
We claim that in this case, 
\begin{equation}
\label{5.1}
\begin{aligned}
      & n\log\rho-(n-1)\log\rho\geq \mu_1(n\log\rho)-\mu_2((n-1)\log\rho)\\
    & \geq a(1-\varepsilon_1)
   (\frac1{2\pi}) (\log|f(\rho^n)-\log|f(\rho^{n-1})| )
    -4a.
\end{aligned}
\end{equation}

In fact, since $\rho>e^{5\pi}$, if  (\ref{Ahlfors1}) does not hold, then (\ref{5.1}) holds
automatically.  Otherwise (\ref{5.1}) follows from (\ref{Ahlfors4}).
We rewrite (\ref{5.1})
$$
\log|f(\rho^n)-\log |f(\rho^{n-1})|<\frac{2\pi}{a(1-\varepsilon_1)}\log\rho
+\frac{8\pi}{(1-\varepsilon_1)}.
$$
so that
\begin{equation}
\label{5.2}
|f(\rho^n)|<C_2\rho^{C_1}|f(\rho^{n-1})|,
\end{equation}
where $C_1$ is as in (\ref{4.0}) and $C_2$ in (\ref{4.2}).
From (\ref{5.2}) it follows that
\begin{equation}
\label{5.2a}
\begin{aligned}
&\frac{|f(\rho^n)|}{\rho^{n(C_1+\varepsilon_2)}}-\frac{|f(\rho^{n-1})|}{\rho^{(n-1)(C_1+\varepsilon_2)}}
<\frac{|f(\rho^n)|}{\rho^{n(C_1+\varepsilon_2)}}-\frac{|f(\rho^{n-1})|}{\rho^{n(C_1+\varepsilon_2)}}\\
&(C_2\rho^{C_1}-1)\frac{|f(\rho^{n-1})|}{\rho^{n(C_1+\varepsilon_2)}}
<\frac{C_2\rho^{C_1}}{\rho^{C_1+\varepsilon_2}}\frac{|f(\rho^{n-1})|}{\rho^{(n-1)(C_1+\varepsilon_2)}}
\end{aligned}
\end{equation}

{\bf Case 2.  } If, on the other hand, the maximum of $\log|h'|$ in $S_n$ is less than $M_0$,
\begin{equation}
\label{5.3}
|f(\rho^n)|-|f(\rho^{n-1})|\leq \int_{\rho^{n-1}}^{\rho^n}(|h'(t)|+1/|h'(t)|)dt
\leq 2e^{M_0}(\rho^n-\rho^{n-1}).
\end{equation}
From (\ref{5.3}) we may write
\begin{equation}
\label{5.4}
\frac{|f(\rho^n))|}{\rho^{n(C_1+\varepsilon_2)}}-\frac{|f(\rho^{n-1)})|}{\rho^{(n-1)(C_1+\varepsilon_2)}}\\
<\frac{|f(\rho^n)|}{\rho^{n(C_1+\varepsilon_2)}}-\frac{|f(\rho^{n-1)})|}{\rho^{n(C_1+\varepsilon_2)}}
\leq 2e^M\frac{(\rho^n-\rho^{n-1})}{\rho^{n(C_1+\varepsilon_2)}}
\end{equation}

Combining (\ref{5.2a}) and (\ref{5.4})  we may write
\begin{equation}
\label{5.4a}
\frac{|f(\rho^n)|}{\rho^{n(C_1+\varepsilon_2)}}
-\frac{|f(\rho^{n-1})|}{\rho^{(n-1)(C_1+\varepsilon_2)}}
\left(\frac{C_2\rho^{C_1}}{\rho^{C_1+\varepsilon_2}}
+\frac 1{\rho^{C_1+\varepsilon_2}}\right)
<\frac{2e^M(1-1/\rho)}{\rho^{n(C_1+\varepsilon_2-1)}}.
\end{equation}
Using (\ref{4.3}) and summing (\ref{5.4a}) up to $n=N-1$  we have
\begin{equation}
\label{5.6}
\frac{|f(\rho^{N-1}))|}{\rho^{(N-1)(C_1+\varepsilon_2)}}
<|f(1)|+K.
\end{equation}
where $K=\sum_{n=1}^\infty 2e^M(1-\rho)/\rho^{n(C_1+\varepsilon_2-1)}$.

To estimate $|f(\sigma)| $ for $\rho^{N-1}<\sigma<\rho^N$, we first consider 
Case 1. If the condition in (\ref{Ahlfors1}) fails, then
\begin{equation}
\label{5.9}
\log|f(\sigma)|-\log|f(\rho^{N-1})|\leq 4\pi\frac{1+\delta_0}{1-\delta_0}.
\end{equation}
which  with (\ref{5.6}) gives
\begin{equation}
\label{5.10}
\log |f(\sigma)|<\log(|f(1)|+K)+(N-1)(C_1+\varepsilon_2)\log\rho +\frac{4\pi}{1-\varepsilon_1}.
\end{equation}
Otherwise, as in (\ref{5.3})-(\ref{5.6})  we have
\begin{equation}
\label{5.11a}
\frac{|f(\sigma))|}{\rho^{N(C_1+\varepsilon_2)}}<|f(1)|+K.
\end{equation}
Comparing (\ref{5.10}) and (\ref{5.11a}) as before, we deduce that (\ref{5.11a}) holds.

Using the fact that $U(\sigma)= 2\sigma$ (recall (\ref{rescale})), $\sigma>\rho^{N-1}$,
 and $C_1+\varepsilon_2=2+2\varepsilon_2$ we have
$$
\frac{\log|f(\sigma)|}{\log U(\sigma)}<\frac{\log|f(\sigma)|}{\log (2\rho^{N-1})}
<\frac{\log\rho^{N(2+2\varepsilon_2)}}{\log (2\rho^{N-1})}+
\frac{\log (|f(1)|+K)}
{\log (2\rho^{N-1})}.
$$
Since $\varepsilon_2$ can be made arbitrarily small, (\ref{mingrowth}) follows.  \qed

\bibliographystyle{amsplain}

\end{document}